
\input vanilla.sty
\magnification 1200
\baselineskip 18pt
\input definiti.tex
\input mathchar.tex
\define\pmf{\par\medpagebreak\flushpar}

\define\e{\varepsilon}

\define\pbf{\par\bigpagebreak\flushpar}

\pmf
\title
Geometry of Cyclic Quotients, I: \\
Knotted Totally Geodesic Submanifolds in Positively Curved Spheres
\endtitle
\author
Alexander Reznikov
\endauthor
\centerline{\bf July, 1994}

\proclaim{\bf Theorem 0.1}  Given a pair of coprime integers
$(m, n), | m |, |n| \ge 2$,
there exists a metric of positive sectional curvature in a compact
three-manifold $M$, with
the following properties:

(a) $M$ is diffeomorphic to $S^3$

(b) There exist geodesic embeddings of $S^1$ in $M$, isotopic to  torus
knots $(m, n)$ and $(m - n)$.
\endproclaim

In connection to the Theorem 0.1 above, we would like to propose the following
problem.
\pmf
{\bf 0.2.} \  Is it true, that a positively curved metric in $S^3$ admits only
finitely many different geodesic knot types?

The problem 0.2 seems to be conceptually related to a theorem of Choi-Schoen on
compactness of the space of embedded minimal surfaces of a given
genius in a Ricci-positive three-manifold (see 3.4.)

\subheading{1. Topology of cyclic quotients}

{\bf 1.1.} \ Here we collect the facts we need on the topology of cyclic
quotients.Let $N$ be a smooth manifold with a smooth action of a cyclic
group $\bbz_n$.  Assume that

(a)  All stationary subgroups of points in $N$ are either trivial or
$\bbz_n$ itself (this is automatically so, if $n$ is prime)

(b) All components of the fixed point set  Fix $(N)$ are of codimension 2.

Then the quotient $N / \bbz_n$ has a canonical structure of a smooth
manifold.  Indeed, let $Q \subset Fix (N)$ be a connected component.
Fix an invariant Riemannian metric in $N$.  Let $\eta$ be the rank
two normal bundle to $Q$ and let $\e D$ be the disc bundle $\eta$ of radius
$\e$.
For $\e$ small the exponential map establishes a diffeomorphism of $\e D$
onto a tubular neighbourhood of $Q$.  We may assume $\eta$ to be orientable
(see below).
Consider $\eta$ as a complex line bundle.  Then cut Exp$(\e D)$ off and
glue the unit
disc bundle of $\eta^{\otimes n}$ to $(N \setminus \text{Exp} (\e D)) /
\bbz_n$.  Doing
this simultaneously for all $Q$, we get a new manifold, homeomorphic to
$N / \bbz_n$.  If $\eta$ is not orientable, consider the double
covering $\tilde Q \overset
\pi \to \rightarrow Q$ such that $\tilde \eta = \pi^\ast \eta$ is orientable.
Denote
$\tau$ the involution of $\tilde Q$ and $\tilde \eta$ such that
$\tilde Q / \tau = Q$ and $\tilde \eta / \tau = \eta$.  Observe that $\tau$
is an orthogonal and antilinear automorphism of $\tilde \eta$, so it induces
an orthogonal and antilinear automorphism of $\tilde \eta^{\otimes n}$
which we denote again by $\tau$.  Now, glue $(\tilde \eta^{\otimes n}) / \tau$
to $(N \setminus \text{Exp} (\e D)) / \bbz_n$.

\demo{1.2 Example}  Let $N$ be a smooth quasiprojective variety over $\bbr$.
Let
$\tau$ be the canonical involution of $N (\bbc)$, coming from
Gal $(\bbc / \bbr)$.
Then $N (\bbc)/ \tau$ is a (real) smooth manifold.  In particular,
$\bbc P^2 / \tau$ is a four-sphere, [5], [6].

\demo{1.3. Cyclic quotients of the three -sphere}  Consider $S^3$ as a unit
sphere in the Hermitian space $\bbc^2$ with the coordinates $(z_1, z_2)$.
Denote
by $K$ and $L$ the geodesic circles $z_2 = 0$ and $z_1 = 0$. For $a, b > 0$
with $a^2 + b^2 = 1$ consider
a torus $T_{a, b} : | z_1 | = a, | z_2| = b$.  The family $T_{a, b}$ form a
fibration of $S^3 \setminus (K \cup L)$, ``converging'' to $K$ and $L$.
Now, any Hopf circle, i.e. an intersection of $S^3$ with a complex
line, lies in one of $T_{a,b}$, namely, $\{ z_2 = \lambda z_1 \}
\cap S^3$ lies on $T_{a,b}$ with $a = \frac{1}{\sqrt{ 1 + | \lambda|^2}}$ and
$b = \frac{\lambda}{\sqrt{1 + | \lambda|^2}}$.  Observe that all tori
$T_{a, b}$ are equidistant from $K$ and $L$ in spherical metric.  Any
such torus is a Heegard surface of the decomposition $S^3 = D^2 \times S^1 \cup
S^1 \times D^2$.

Now, consider a $\bbz_m$-action $(z_1, z_2) \mapsto (z_1, e^{\frac{2 \pi k}{m}}
z_2)$.  It
has $K$ as a fixed point set and acts free in $S^3 \setminus K$.  According
to 1.1, $S^3 / \bbz_m$ is a manifold.  We claim $S^3 / \bbz_m \approx S^3$.
Indeed,
the action in the handle, which contains $L$, is free and the quotient is
obviously a handle again.  The action in the other handle, which contains $K$
is
fiber like, as in 1.1, and since the normal bundle to $K$ is trivial, the
quotient is again a handle, which proves the statement.

Let $n$ be an integer, coprime to $m$.  Consider action of $\bbz_m
\times \bbz_n$ by $(z_1, z_2) \mapsto
(e^{\frac{2 \pi i k}{m}} z_1, e^{\frac{2 \pi k}{n}} z_2)$.  Applying the
diffeomorphism above
twice, we get a following lemmas.

\proclaim{Lemma 1.4}  The quotient $S^3 / \bbz_m \times \bbz_n$ is
diffeomorphic
to $S^3$.
\endproclaim

\subheading{2. Constructing a metric in the cyclic quotient}

The construction below uses the computations of Gronov and Thurston [4].
In that paper, Gromov and Thurston introduced negatively curved metrics
on ramified
coverings of hyperbolic manifolds.  Our situation is ``dual'' to that
considered in [4],
in particular, lifting to a ramified covering is replaced by
descending to a cyclic
quotient.

\proclaim{Lemma 2.1}  (comp [4], p.4).  Given $n \in \bbn$ and
$\rho > 0$, there
exists a smooth function $\sigma (r)$ with the following properties:

(i) $\sigma (r) = \sin r$ for small $r$

(ii) $\sigma' (r) > 0$ and $\sigma^{''} (r) < 0$

(iii) $\sigma (r) = \frac{\sin r}{n}$ for $r \ge \rho$
\endproclaim

The proof is immediate.

Now, the metric of $S^3 \setminus L$ can be written as
$$ g = d r^2 + \sin^2 r d \theta^2 + \cos^2 r \ d t^2 $$
Here $t$ is the length parameter along $K$ and $(r, \theta)$ are
polar coordinates
in geodesic two-spheres, orthogonal to $K$.

Consider the cyclic quotient $S^3 / \bbz_n$ and equip it with the metric
$$ \tilde g = d r^2 + \sigma^2 (r) d \theta^2 + \cos^2 r d t^2 $$
where $\theta$ here is the \underbar{new} angle parameter.  Outside the
$\rho$-neighbourhood of $K /\bbz_n$, this
is just a descend of the spherical metric by the (free) action of $\bbz_n$.
The crucial fact is the following.

\proclaim{Lemma 2.2}  The metric $\tilde g$ is a well-defined smooth
metric on $S^3 / \bbz_n$, of strictly positive
curvature, which is invariant under the (descend of) the
$\bbz_m$-action.  Out of the small
neighbourhood of $K$, the curvature of $\tilde g$ is constant.
\endproclaim

\demo{Proof}  It is elementary to check that $\tilde g$ is smooth with
respect
to the manifold structure of $S^3 / \bbz_n$.  The positivity of curvature
follows form computations of [4], p. 4-5, with obvious changes
($\cos h \to \cos$ etc.).
the invariance under the $\bbz_m$-action is obvious from the construction.

\demo{2.3}  Taking $\rho$ small and repeating the construction with
respect to
the $\bbz_m$-action, we come to the following lemma.

\proclaim{Lemma 2.3}  The quotient $S^3 / \bbz_m \times \bbz_n$ can
be equipped
with the metric $\overset \approx \to g$ with following properties:

(a) the curvature of $(S^3 / \bbz_m \times \bbz_n, \overset
\approx \to g)$ is
strictly positive

(b) outside arbitrary small neighbourhood of $K / \bbz_m$ and $L /
\bbz_n$ the metric $\overset
\approx \to g$ is a descend of the spherical metric of $S^3$.
\endproclaim

\subheading{3. Knotted geodesics}

\proclaim{Lemma 3.1}  The image of any Hopf geodesic circle in
$S^3$ outside the $\rho$-neighbourhoods of $K, L$ is a torus
knot $(m, n)$ in $S^3 / \bbz_m \times \bbz_n \approx S^3$.
\endproclaim

\demo{Proof}  Let $\gamma = (z_2 = \lambda z_1) \cap S^3, \lambda
\in \bbc$, be a Hopf circle.
According to 1.3, $\gamma \subset T_{a,b}$ with $a =
\frac1{\sqrt{1 + | \lambda|^2}}$.
In angle coordinates $(e^{i \theta}, e^{i \tau})$ on $T_{a,b}, \ \gamma$
may be written in parametric
form as $t \mapsto (e^{i t}, e^{i t})$.  The quotient map
$T_{a,b} \to T_{a,b} / \bbz_m \times \bbz_n \approx T^2$ can
be written as $(e^{i \theta}, e^{i \theta}) \mapsto
(e^{m i \theta}, e^{n i \tau})$.  So the image of $\gamma$ is
$t \mapsto (e^{m i t}, e^{n i t})$ which is $(m, n)$ torus knot.

\demo{3.2. Geodesics of different knot type}  Consider a new complex
structure in $\bbc^2$, defined by the matrix $\pmatrix -i &0 \\
0 &i \endpmatrix $.  Observe that $K$ and $L$ are still Hopf circles which
respect to this new complex structure.  The tori $T_{a,b}$ have
the same equation
$| w_1 | = a, \ | w_2 | = b$ with respect to new coordinates
$w_1 = z_1, \ w_2 = \bar z_2$.
Hence they contain a Hopf circle $w_2 = \lambda w_1$, which descends
to a geodesic in $M$, which is isotopic to
the torus knot $(m, - n)$.  Since torus knots are
invertible, but not amphicherical ([3]), we get two different
knot types among geodesics in $M$.  This concludes the proof of
the Theorem 0.1.

\demo{3.3. Questions and Remarks 3.3.1.}  May knots other
than torus knots be realized as
geodesics of positively curved metric in $S^3$?

\demo{3.3.2.}  Fix a pinching constant $\delta$.  Can a three-sphere with
a $\delta$-pinched metric have geodesics of arbitrary torus knot type?

\demo{3.4.} Suppose $S^g \subset (S^3$, can) be   an immersed
minimal surface of
genus $g$ whose image do not touch $K \cup L$.  Applying the construction
above, we
come to a minimal surface in $M$ with ``a lot of'' selfintersections.
In particular, one may start with a Clifford
torus, close to $T_{1/ \sqrt{2}, 1 \sqrt{2}}$.  This situation contrasts the
compactness theorem of \underbar{embedded} minimal surface of a given genus
[7].  For
$g \ge 2$, are may therefore ask a following question:

May a compact minimal surface in $S^3$ of genus $g \ge 2$ avoid a
geodesic circle?
\pmf
{3.5.}  We sketch a different type of examples which lead to
Theorem 0.1.  in case when both $m, n$ are odd.  Start
with a positively curved metric in $S^2$.  Let $\bar M =
U S^2$, a unit tangent bundle with the Sasaki metric
(making $U S^2 \to S^2$ be
a Riemannian submersion).  If the curvature of $S^2$ is less
than $1 / \sqrt{3}$,
then $U S^2$ is positively curved.  This follows from the
O'Neil formulas for
Riemannian submersions with totally geodesic fibers
([2],chap.9).   Since
$U S^2 \cong \bbr P^2$, the double cover of $U S^2$ is the three-sphere.

Now, we may take the metric of $S^2$ to be  rotationally invariant.  Then
the
computations of ([1]) show that we have the full
control on closed geodesics of $S^2$.  In
particular, there are closed trajectories of the geodesic flow in
$U S^2$, whose lift to $S^3$ will
be a torus knot of a given type $(m, n)$ if both $m, n$ are odd.
Unfortunately,
already the trefoil knot may not be realized in this way.

\demo{3.6. Concluding remarks}  It looks like that there exists a
positively curved metric on $S^4$ with a totally geodesic $\bbr P^2$
having the
normal Euler number four.  Indeed, look at the standard K\"ahler metric
in $\bbc P^2$.  The canonical
autoholomorphic
involution $\tau: \bbc P^2 \to \bbc P^2$ is an isometry and
$\bbc P^2 / \tau \cong S^4$.  The
fixed point set Fix $(\tau)$ is a totally geodesic $\bbr P^2
\subset \bbc P^2$.  It is
possible to mimic the construction of 2.2. and find a perturbation
of the quotient metric in
$\bbc P^2 / \tau$ in directions, orthogonal to $\bbr P^2$ (recall
that there exists exactly one totally geodesic surface, also
isometric to $\bbr P^2$, meeting Fix $(\tau)$ orthogonally at
a given point.)  The curvature tensor, however, is no more diagonal , and
it is a nontrivial problem to check if the curvature is positive.
Observe that the resulting metric admit an isometric $SO (3)$-action,
and the equidistant from $\bbr P^2$ manifolds are lense spaces $S O (3) /
\bbz_2 \cong S^3 / \bbz_4$ with
a homogeneous metric, which is different from Berger's metrics.
\pbf
\centerline{\bf References}

\item{1.} A. Besse, {\it Manifolds, All of Whose Geodesics are Closed},
Springer, 1978.

\item{2.} A. Besse, {\it Einstein Manifolds}, Springer, 1987.

\item{3.} G. Burde, H. Zieschang, {\it Knots}, Walter de Gruyter, 1985.

\item{4.} M. Gromov, W. Thurston, {\it Pinching constants for hyperbolic
manifolds},
Invent. Math. {\bf 89} (1987), 1--12.

\item{5.} M. Kreck, {\it On the homeomorphism classification of smooth
knotted surfaces in the 4-sphere} in: S. K. Donaldson and C. B. Thomas, ed.,
{\it Geometry of
Low-dimensional Manifolds}, I, Cambridge Univ. Press. 1990.

\item{6.} L. Guillou, A.Marin, {\it A la Recherche de la Topologie Perdue},
Birkh\"auser, 1986.

\item{7.} H.I.Choi, R.Schoen, {\it The space of minimal embeddings of a surface
into a three-dimensional manifold of positive Ricci curvature}, Inv. Math.{\bf
81} (1985), 387--394.
\pbf
Department of Mathematics
\par
\noindent Hebrew University
\par
\noindent Givat Ram 91904
\par
\noindent Israel
\par
\noindent email: simplex\@math.huji.ac.il
\bye